\documentclass[10pt]{article}

\usepackage{latexsym}
\usepackage{enumerate, amssymb, amsmath, amsthm}
\usepackage[utf8]{inputenc}
\usepackage{graphicx}
\usepackage{url}
\usepackage{color}

\theoremstyle{theorem}

\theoremstyle{definition}

\title{A probabilistic way to discover the rainbow}
\author{Joscha Prochno and Michael Schmitz}
\date{}

\begin{document}
\maketitle

\section{Introduction}
A slogan that you find on the back of a pack of Skittles candy says ``No two rainbows are the same. Neither are two packs of Skittles. Enjoy an odd mix.''. In the online blog \cite{PossiblyWrong} it is described how the blogger found two identical packs of Skittles, among 468 packs with a total of 27,740 Skittles. Meticulously collecting the data for this experiment was apparently triggered by some earlier calculations. More precisely, the blogger writes: 

\begin{center}
	``A few months ago, we did some calculations on a cocktail napkin, so to speak, predicting that we should be able to find a pair of identical packs of Skittles with a reasonably -- and perhaps surprisingly -- small amount of effort.''
\end{center} 

Whether performing this admittedly yummy experiment really only requires a ``small amount of effort'' or not, we chose to write this short article as intellectual candy for the mathematically inclined reader and model this experiment as a probabilistic one. This allows to quantify the probability that two randomly selected packs of Skittles candy are identical and, in a next step, to estimate the expected number of packs one has to purchase until the first match. The approach requires merely elementary probability theory and, as is typical for such a discrete problem, some combinatorial considerations. 

We also believe this problem to be appealing for middle and high school students as such an experiment can be repeated and subsequently be analysed by probabilistic tools.  Here, an adaption in the precision of arguments creates a certain variability in the level of requirements. Concluding the article more sophisticated tools, such as generating functions, are employed, which certainly exceed middle or high school level but are suitable for working with undergraduate university students.

In order to approach this question mathematically, we need to start with a suitable model. Let us assume that each pack of Skittles contains the exact same number $n\in\mathbb{N}$ of Skittles (for brevity we sometimes say $n$-packs of Skittles) and that there are $d\in\mathbb{N}$ different colours.\footnote{Actually the number varies from pack to pack, in \cite{PossiblyWrong} it says that most studies suggest an average of about 60 candies per pack. There are five different colours/flavours.} When filling a pack of Skittles we want to assume that we do it randomly and with a uniform distribution over the $d$ colours, i.e., for each colour $k \in \{1,\ldots, d\}$ and each Skittle entering a pack, the probability that it has the colour $k$ is given by $\frac{1}{d}$. We shall say that two packs of Skittles are identical if for each colour they contain the same number of Skittles.

\section{Two colours}\label{Two colours}

A key to solving the Skittles problem elegantly will be a reinterpretation that can be best understood when considering only $d=2$ colours for the time being, and generalize the idea to an arbitrary number $d\leq n$ of colours afterwards. We consider the filling process of a pack of $n$ Skittles with two possible colours, say red and green, and imagine it to be constructed by a planar random walk from the origin (since we start with an empty pack) on the integer lattice in the plane, where a step right adds a red sweet, and a step up adds a green sweet; each of the two possibilities being equally likely, having probability $1/2$, while the steps are independent. The pack is full on the line $x + y = n$, and the finishing point on that line corresponds to a pack of $x$ reds and $y$ greens. The following picture shows two lattice paths (dashed/solid) that both correspond to a pack of eleven Skittles with six reds and five greens.

\begin{center}
	\includegraphics[width=0.8\textwidth]{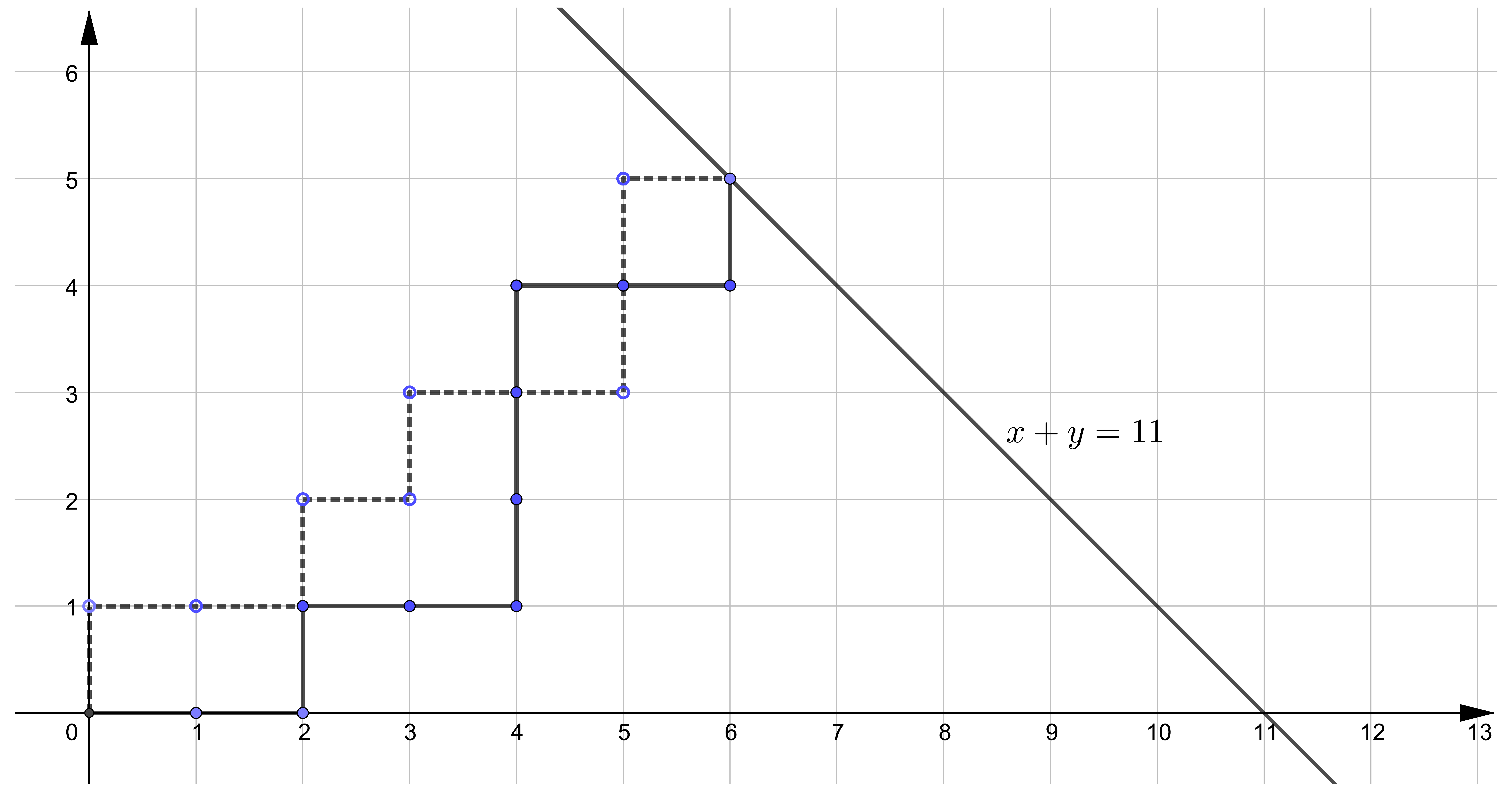}\\
	{\small \textbf{Picture 1:} Lattice paths represent Skittles packs.}
\end{center}

There are $\binom{n}{x}$ $(=\binom{n}{y})$ paths from the origin to a given point $(x,y)$ with $x+y=n$, as we can choose any of the $n$ total steps to be the $x$ steps to the right (or the $y$ steps up). Note that since $x+y=n$, indeed
\[
{n\choose x} = {n\choose n-y} = {n\choose y},
\]
where we used the identity ${n\choose k} = {n\choose n-k}$, which holds for all $n\in\mathbb N$ and $k\in\{0,1,\dots,n\}$.
Clearly, each path with $n$ steps has the probability $(1/2)^n$ and therefore such a random walk ends at $(x,y)$ with probability $\binom{n}{x}/2^n$.

Thus, two independent random walks both end at a given point $(x,y)$ with the probability $\binom{n}{x}^2/2^{2n}$. Summing over all possible $x\in\{0,1,\dots,n\}$ gives the probability of the event $E_n^2$ that two independent random walks end at a mutual point with step distance $n$ from the origin, namely
\begin{equation}\label{prob_of_E_n^2}
	\mathbb{P}[E_n^2] = \frac{\sum_{x=0}^n \binom{n}{x}^2}{2^{2n}}.
\end{equation}
Note that we have hereby solved the problem at hand, because $\mathbb{P}[E_n^2]$ equals the probability that two independently chosen $n$-packs of Skittles with two possible colours are identical, although the cardinality of the event $E_n^2$ is not the number of pairs of identical $n$-packs with two possible colours. The trick is that we chose a model that takes order into account (by considering paths: each step represents a sweet entering a pack), although there is clearly no order in a pack of Skittles.\\

We can use a reinterpretation to simplify (\ref{prob_of_E_n^2}). To this end we consider the same kind of random walks with $2n$ instead of $n$ steps. The probability that such a walk ends at $(n,n)$ equals $\binom{2n}{n}/2^{2n}$, as $n$ steps to the right have to be chosen out of $2n$ steps. On the other hand, such a walk has to pass the line $x+y = n$ at some point $(x,y)$. This means that amongst the first $n$ steps $x$ steps to the right have been made. Then, in order to end at $(n,n)$, amongst the remaining $n$ steps $n-x$ steps to the right have to be made (see picture 2).

\begin{center}
	\includegraphics[width=0.8\textwidth]{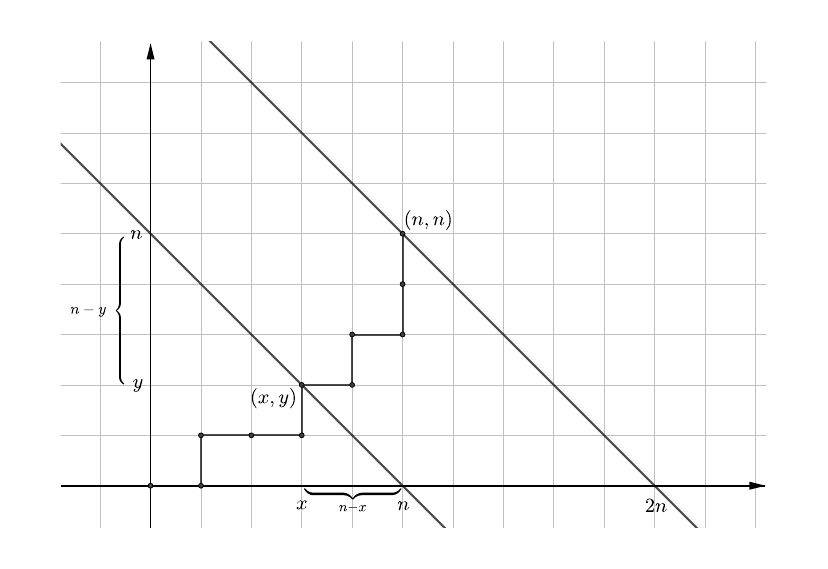}\\
	{\small \textbf{Picture 2:} A walk with $2n$ steps passes a point $(x,y)$ with $x+y=n$.}
\end{center}

Therefore, there are $\binom{n}{x} \cdot \binom{n}{n-x} = \binom{n}{x}^2$ paths that end at $(n,n)$ and pass through the point $(x,y)$. Summing over all $x$ tells us that there are $\sum_{x=0}^n \binom{n}{x}^2$ walks of step length $2n$ that end at $(n,n)$,\footnote{What we did here is essentially the combinatorial proof of the Vandermonde identity ${m_1+m_2\choose n} = \sum_{k=0}^n {m_1\choose k}{m_2\choose n-k}$ for the case $m_1=m_2=n$.} and we obtain
$$\mathbb{P}[E_2^n] = \frac{\sum_{x=0}^{n} \binom{n}{x}^2}{2^{2n}} = \frac{\binom{2n}{n}}{2^{2n}}.$$

\section{Three and more colours}\label{Three and more colours}

Now we consider $n$-packs of Skittles with $d$ possible colours, and start to generalize our approach by tackling the case $d=3$. That is, considering a spatial random walk from the origin on the integer lattice in $3$-dimensional space, where a step right adds a red sweet, a step forward adds a green sweet, and a step up adds a blue sweet. The pack is full on the plane $x+y+z=n$. Analogously to above we denote the event that two independently performed $n$-step random walks end at a mutual point by $E_n^3$ (or $E_n^d$ in general).

\begin{center}
	\includegraphics[width=0.55\textwidth]{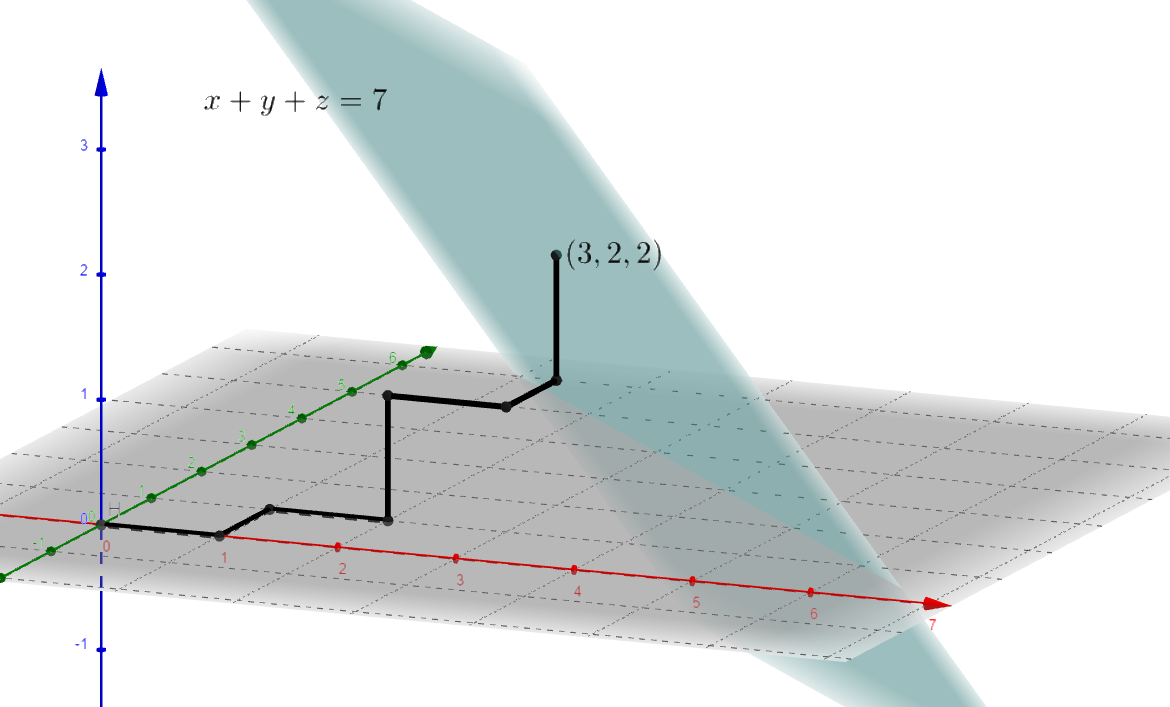} \includegraphics[width=0.44\textwidth]{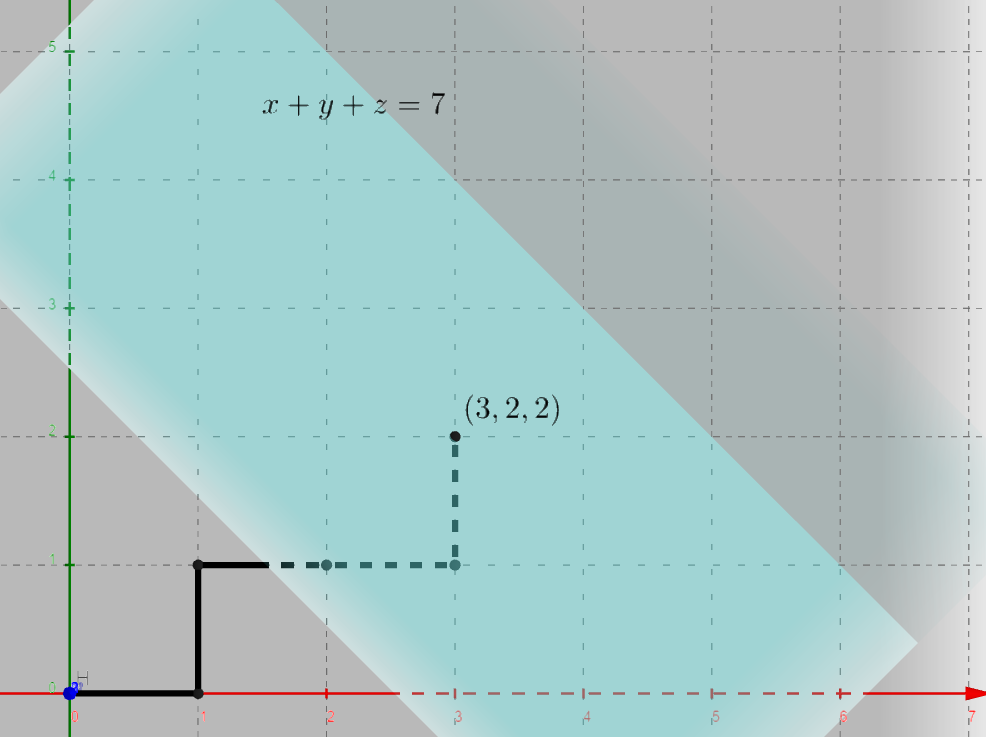}\\
	{\small \textbf{Picture 3:} A $7$-step spatial random walk and its orthogonal projection from above.}
\end{center}

For a walk to end at a given point $(x,y,z)$ with $x+y+z=n$, there have to be made $z$ steps up and we have $x+y=n-z$. This means that if we look at the situation directly from above (i.e., consider an orthogonal projection on the $xy$-plane), we have a random walk with $n-z$ steps in the plane that ends at the point $(x,y)$, so we are back in the situation considered earlier (see picture 3). Thus, the number of pairs of walks that both perform exactly $z$ steps up and have the same endpoint equals $\binom{n}{z}^2 |E_{n-z}^2|$, and summing over all possible $z\in\{0,1,\dots,n\}$ yields
\begin{equation}\label{rec_E_n^3}
	|E_n^3| = \sum_{z=0}^{n} \binom{n}{z}^2 |E_{n-z}^2|.
\end{equation}
As each path has probability $(1/3)^n$, we obtain $\mathbb{P}[E_n^3] = |E_n^3|/3^{2n}$. Formula (\ref{rec_E_n^3}) is a nice recursion, but if we pause for a moment, we see that the same considerations also provide a non-recursive expression. For a walk to end at $(x,y,z)$ we have $\binom{n}{z}$ choices for when to make the $z$ steps up, and out of the remaining $n-z$ steps we have $\binom{n-z}{x}$ choices when to make the $x$ steps to the right. Therefore, there are $\binom{n}{z}^2 \binom{n-z}{x}^2$ pairs of walks with $z$ steps up and $x$ steps to the right (and consequently $y$ steps forward). Summing over all possible $x,y,z$ now yields
$$|E_n^3| = \sum_{x+y+z=n\atop{x,y,z \in \{0,1,\ldots, n\}}}  \binom{n}{z}^2 \binom{n-z}{x}^2.$$
Keeping the assumption $n-z-x=y$ in mind, we observe that
$$\binom{n}{z} \binom{n-z}{x} = \frac{n!}{(n-z)!z!} \cdot \frac{(n-z)!}{(n-z-x)!x!} = \frac{n!}{z! x! y!},$$
and recognize the multinomial coefficient\footnote{The general definition is $\binom{n}{x_1,\ldots, x_d}:= \frac{n!}{x_1!\ldots x_d!}$. A well known and nice interpretation of the multinomial coefficient is an alphabetical jumble, i.e., the number of distinct permutations of a word of length $n$ in which $d$ different letters occur and each letter $i$ occurs $x_i$ times. For instance, there are $\binom{8}{2,2,1,1,1,1}$ different permutations of the Word SKITTLES.} $\binom{n}{x,y,z} := \frac{n!}{x! y! z!}$. Therefore, we may rewrite
$$|E_n^3| = \sum_{x+y+z=n\atop{x,y,z \in \{0,1,\ldots, n\}}} \binom{n}{x,y,z}^2.$$

It is now not so hard to generalize these ideas to $n$-packs of Skittles with an arbitrary number $d$ of possible colours. Filling a pack can be thought of as a spatial random walk from the origin on the integer lattice in $d$-dimensional space, where a step in $x_1$-direction adds a sweet of colour 1, a step in $x_2$-direction adds a sweet of colour 2, and so on. Now, the considerations for both the recursive and for the non-recursive formula are pretty much as above.

For the recursion we fix a number $k$ of steps that are made in direction $x_1$. There are $\binom{n}{k}$ possible choices for these $k$ steps. The remaining $n-k$ steps have to be carried out in $d-1$ dimensions, so there are $\binom{n}{k}^2|E_{n-k}^{d-1}|$ pairs of walks that perform exactly $k$ steps in direction $x_1$ and end at a mutual point. Again, we sum over all $k$ and obtain
\begin{equation}\label{general_recursion}
	|E_n^d| = \sum_{k=0}^{n} \binom{n}{k}^2 |E_{n-k}^{d-1}|.
\end{equation}
The non-recursive formula is also generalized in a straightforward manner. For a random walk to end at a given point $(k_1,k_2,\ldots k_d)$ in $d$-dimensional space with $k_1+k_2+ \ldots + k_d = n$, we have to choose $k_1$ steps in $x_1$-direction, $k_2$ steps in $x_2$-direction, and so on. There are $\binom{n}{k_1}$ possible choices for the $k_1$ steps in $x_1$-direction. From the remaining $n-k_1$ steps we have $\binom{n-k_1}{k_2}$ possibilities to choose $k_2$ steps in $x_2$-direction, and so on. Therefore, there are
$$\binom{n}{k_1} \cdot \binom{n-k_1}{k_2} \cdot \ldots \cdot \binom{n-[k_1+\ldots+k_{d-1}]}{k_d} = \binom{n}{k_1,k_2, \ldots, k_d}$$
paths with $n$ steps that end at $(k_1,k_2, \ldots, k_d)$. Summing over all possible $k_1,\ldots, k_d$ yields
\begin{equation}\label{general_formula}
	|E_n^d| = \sum_{k_1+\dots+k_d=n\atop{k_i\in\{0,\dots,n\}}} {n \choose k_1,\dots,k_d}^2.
\end{equation}

As each path occurs with probability $(1/d)^n$, we obtain $\mathbb{P}[E_n^d] = |E_n^d|/d^{2n}$.

\section{Crunching some numbers}
Let us use our formulas to compute some values and see why it is nice to have both closed and recursive expressions. By (\ref{general_formula}) we have
$$|E_2^2| = \sum_{k_1+k_2=2\atop{k_i \in \{0,1,2\}}} \binom{2}{k_1,k_2}^2 = \binom{2}{0,2}^2 + \binom{2}{1,1}^2 + \binom{2}{2,0}^2 = 6.$$
That was easy, so we could try dealing with slightly larger numbers in another example, e.g.,
$$|E_3^3| = \sum_{k_1+k_2+k_3=3\atop{k_i \in \{0,1,2,3\}}} \binom{3}{k_1,k_2,k_3}^2.$$
We can arrange the sum $3=0+0+3$ in three and $3=0+1+2$ in six possible orders, while $3=1+1+1$ has only one possible order. Therefore, we obtain
$$|E_3^3| = 3 \binom{3}{0,0,3}^2 + 6\binom{3}{0,1,2}^2 + \binom{3}{1,1,1}^2 = 3+6\cdot 9 + 36 =93.$$
Considering this (for $n=3$ and $d=3$ the sum already consists of ten summands!) we are lucky to have a recursion for determining (e.g., by means of a computer) the numbers $|E_d^n|$ for larger $n$ and $d$. Tables 1 and 2 show the values for $|E_d^n|$ and $\mathbb{P}[E_d^n]$ (rounded to four digits) for $1\le n,d \le 5$.

\begin{center}

	\begin{tabular}{|c||c|c|c|c|c|}
		\hline
		n \textbackslash d & 1 & 2 & 3 & 4 & 5 \\
		\hline\hline
		1 & 1 & 2 & 3 & 4 & 5 \\
		\hline
		2 & 1 & 6 & 15 & 28 & 45 \\
		\hline
		3 & 1 & 20 & 93 & 256 & 545 \\
		\hline
		4 & 1 & 70 & 639 & 2716 & 7885 \\
		\hline
		5 & 1 & 252 & 4653 & 31504 & 127905 \\
		\hline
	\end{tabular}\\
	{\small Table 1: $|E_d^n|$ for $1\le n,d \le 5$}\\[4mm]

	\begin{tabular}{|c||c|c|c|c|c|}
		\hline
		n \textbackslash d & 1 & 2 & 3 & 4 & 5 \\
		\hline\hline
		1 & 1 & 0.5 & 0.3333 & 0.25 & 0.2 \\
		\hline
		2 & 1 & 0.375 & 0.1825 & 0.1094 & 0.072 \\
		\hline
		3 & 1 & 0.3125 & 0.1276 & 0.0625 & 0.0349 \\
		\hline
		4 & 1 & 0.2734 & 0.0974 & 0.0414 & 0.0202 \\
		\hline
		5 & 1 & 0.2461 & 0.0788 & 0.0300 & 0.0131 \\
		\hline
	\end{tabular}\\
	{\small Table 2: $\mathbb{P}[E_d^n]$ for $1\le n,d \le 5$}

\end{center}

Of course, we also want to know the probability that two randomly purchased packs of Skittles are identical assuming the realistic values of $d=5$ colours and $n=60$ sweets in each pack. This is $\mathbb{P}[E_5^{60}] = 0.00009752... \approx 0,01\%$. So while it is obviously not true that no two packs of Skittles are same (as the slogan claims), it is at least very unlikely (in particular, because not all packs contain the same number of Skittles, which decreases the probability even further!).

\section{Expected number of packs needed for a match}\label{sec:expected_number}

We now imagine that somebody purchases a pack of Skittles candy each day and compares it to any of the previously bought packs to see if it is identical to one of them (as said above, there are actually people who do such things). We now ask the following question:
\begin{center}
	\emph{How many packs must be bought on average until the first match appears, i.e., what is the expected value of purchased packs in this experiment?}
\end{center}

First, we present a plausible, yet faulty approach, and then discuss where the error happened, why the question is not so easy to answer precisely, and how to get at least a good estimate for the desired expected value. We know that two independent $n$-step random walks $W, W'$ represent identical $n$-packs of Skittles if and only if they have the same endpoint. In this case let us say that they are \textbf{equivalent} and write $W \sim W'$. In order to tackle the problem at hand, we switch from considering two walks to considering a sequence $W_1, W_2, W_3, \ldots$ of independent $n$-step random walks. In other words, $W$ is the $n$-step random walk in $\mathbb{R}^d$ that ends at some point on the hyperplane $H=\{(x_1,\ldots, x_d)\in \mathbb{N}_0^d~|~ x_1+\ldots+x_d=n\}$ considered above, and $W_1,W_2,\ldots$ are independent, identically distributed copies of $W$. As we are only interested in whether the random walks have the same endpoint or not, we can switch to considering random variables $X_1,X_2,\ldots$ with values in $H$, where $X_i=(k_1,\ldots,k_d)$ if and only if $W_i$ ends in $(k_1,\ldots,k_d)$. From our considerations above, we know that $\mathbb{P}[X_i=(k_1,\ldots, k_d)]=\frac{\binom{n}{k_1,\ldots,k_d}}{d^n}$ for every $i \in \mathbb{N}$.

Moreover, we consider a random variable $Z$ that shall return the number of purchases until the first match, i.e., $Z$ takes the value $\ell \ge 2$ if and only if $W_1, W_2, \ldots, W_{\ell-1}$ are pairwise non-equivalent and $W_\ell \sim W_i$ for some $i \in \{1,\ldots, \ell-1\}$. In other words, $\mathbb{P}[Z=\ell]$ is given by
$$\mathbb{P}[X_1, \ldots, X_{\ell-1} \text{ p.d.}] \cdot \mathbb{P}[X_\ell = X_i \text{ for some } i<\ell~|~ X_1, \ldots, X_{\ell-1} \text{ p.d.}],$$
where we abbreviated `pairwise different' by `p.d.' and the latter factor denotes a conditional probability. Now comes the plausible, yet incorrect argument: for $i\ne j$ let $p=\mathbb{P}[X_i=X_j]=\mathbb{P}[E_n^d]$ be the probability that two arbitrary random walks have the same endpoint. In other words, $p$ is the probability that two randomly purchased $n$-packs of Skittles with $d$ possible flavors are identical. Among $\ell-1$ packs of Skittles (or random walks) there are $\binom{\ell-1}{2}$ different pairs, and each pair consists of different packs with probability $1-p$. Thus,
$$\mathbb{P}[X_1, \ldots, X_{\ell-1} \text{ p.d.}] = (1-p)^{\binom{\ell-1}{2}}.$$
Moreover, given that $X_1, \ldots, X_{\ell-1}$ are pairwise different, the probability that $X_{\ell}$ matches with at least one of the previous $X_i,~i<\ell$, equals $(\ell-1)p$. Hence,
$$\mathbb{P}[Z=\ell] = (1-p)^{\binom{\ell-1}{2}} \cdot (\ell-1)p.$$
What did we do wrong? One point is that
$$\mathbb{P}[X_\ell = X_i \text{ for some } i<\ell~|~ X_1, \ldots, X_{\ell-1} \text{ p.d.}] \le 1,$$
but $(\ell-1)p$ clearly does exceed $1$ for large $\ell$. This could be easily taken into account by replacing the factor $(\ell-1)p$ by $\min\{1,(\ell-1)p\}$, but we made another, more critical mistake: the events of having identical packs of Skittles are clearly not independent. For example, the probability that $X_3$ is identical to $X_1$ or $X_2$ depends on whether $X_1$ and $X_2$ are identical or not. One might think that $(1-p)^{\binom{k}{2}}$ should at least serve as an upper bound for $\mathbb{P}[X_1,\ldots,X_k \text{ p.d.}]$, but even this turns out to be not true in general. We want to explain more precisely why the suggested upper bound seems plausible at first and what we mean by `not true in general'. We are considering a series of i.i.d copies of $X:\Omega \to H$. Let us denote $q_x=\mathbb{P}[X=x]$ for each $x=(x_1,\ldots,x_d) \in H$. Then, we have $q_x \in [0,1]$ for each $x$ and $\sum_{x\in H} q_x =1$. If $X$ was uniformly distributed, i.e., $q_x=1/m$ for each $x$, where we let $m=|H|$, we would have
$$\mathbb{P}[X_1, \ldots, X_k \text{ p.d.}] = \frac{\prod_{i=0}^{k-1} (m-i)}{m^{k}}.$$
Note that for $m=365$ this is nothing but the classical birthday problem (see, e.g. \cite{Feller}). Here, it holds $p=1/m$, and it can be shown without much effort that $(1-1/m)^{\binom{k}{2}}$ is indeed an upper bound and a good estimate for the probability that $X_1,\ldots,X_k$ are pairwise different (see, e.g., \cite{Schwarz, Perkins}). For example, with $m=365$ and $k=23$, it holds
$$\frac{\prod_{i=0}^{k-1} (m-i)}{m^{k}} = 0.4927\ldots \quad \text{ and } \quad \left(1-\frac{1}{m}\right)^{\binom{k}{2}}=0.4995\ldots.\footnote{This is the famous fact that among 23 people the chances that at least two of them have the same birthday (day and month) are more than even (assumed that birthdays are uniformly distributed throughout the year).}$$
The fact that for a uniform distribution the faulty, dependence-ignoring approach yields an upper bound and a good approximation makes it plausible that this could also be true for arbitrary distributions, but the following counterexample\footnote{This counterexample was delivered by Will Perkins and pointed out on \cite{mathoverflow}.} shows that this not the case. Let $X:\Omega \to \{1,2,3\}$ with $\mathbb{P}[X=1]=0.8$ and $\mathbb{P}[X=2]=\mathbb{P}[X=3]=0.1$. Then the probability that three copies of $X$ are pairwise different equals $6\cdot 0.8\cdot 0.1^2 =0.048$. On the other hand, the probability $p$ that two arbitrary copies of $X$ coincide is given by $p=0.8^2+0.1^2+0.1^2$, and therefore the would-be upper bound equals $(1-p)^{\binom{3}{2}}=0.039\ldots$.

Another interesting fact is that the probability that $k$ copies of $X$ are pairwise different is maximized by the uniform distribution. This is not so hard to show (see, e.g., \cite{Munford, Nunnikhoven}) and, at first glance, seems promising to allow a simple approach to get an upper bound for the desired expected value. But on second thought, this idea is futile, as $\mathbb{P}[Z=\ell]$ is given by
$$\mathbb{P}[X_1, \ldots, X_{\ell-1} \text{ p.d.}] \cdot \mathbb{P}[X_\ell = X_i \text{ for some } i<\ell~|~ X_1, \ldots, X_{\ell-1} \text{ p.d.}],$$
and the uniform distribution yields an upper bound for the first factor, but a lower bound for the second factor. Fortunately, there are other techniques that solve the problem. In \cite{Wiener} it is shown that
$$\mathbb{P}[X_1, \ldots, X_{\ell-1} \text{ p.d.}] \le \left(1-\sqrt{p}\right)^{k-1}\left(1+\sqrt{p}(k-1)\right),$$
and
$$\sqrt{\frac{\pi}{2p}-\frac{2}{5}} < \mathbb{E}[Z] \le \frac{2}{\sqrt{p}}.$$
For the proofs of these inequalities we refer the interested reader to \cite{Wiener}. Here, we are content to compare the expected value that would result from our intuitively reasonable, yet faulty first approach with a correct estimate for $n=60, d=5$. From $P[Z=\ell] =  (1-p)^{\binom{\ell-1}{2}} \cdot \min\{1,(\ell-1)p\}$ it would follow that $\mathbb{E}[Z] \approx 129$, and Wiener's result in \cite{Wiener} yields
$$ 126.9... < \mathbb{E}[Z] \le 202.5...$$ 
We see that the incorrect approach yields a plausible result. It would be interesting to examine whether this approach can be adapted suitable to obtain a correct estimate in our case, that is, the case of a series of multinomial distributed random variables. Note that the correct upper bound derived from Wiener's results is much smaller than the result of the experiment in the blog \cite{PossiblyWrong}. And in fact it should be, considering	that actually the packs of Skittles may contain different numbers of candies, which clearly reduces the probability of a match and therefore increases the
expected value.

\section{Comparison to other results}

In the internet blog mentioned above a generating function for the numbers $\mathbb{P}[E_d^n]$ is presented. In particular, it says that
$$\mathbb{P}[E_d^n] = \frac{1}{d^{2n}} \left[\frac{x^{2n}}{(n!)^2} \right] \left(\sum_{k\ge 0} \left(\frac{x^k}{k!}\right)^2\right)^d, \qquad (d,n \in \mathbb{N}, d\le n)$$
which means that $|E_d^n|=d^{2n} \mathbb{P}[E_d^n]$ is the coefficient of $\frac{x^{2n}}{(n!)^2}$ in the $d^{\text{th}}$ power of the series $\sum_{k\ge 0} \left(\frac{x^k}{k!}\right)^2$. We want to see that this corresponds exactly to our result. We recall from \eqref{general_formula} that
$$
  |E_d^n| = \sum_{k_1+\dots+k_d=n\atop{k_i\in\{0,\dots,n\}}}\frac{(n!)^2}{(k_1!k_2! \cdots k_d!)^2}.
$$
Thus, it remains to show that the coefficient of $x^{2n}$ in the $d^{\text{th}}$ power of the series $\sum_{k\ge 0} \left(\frac{x^k}{k!}\right)^2$ equals
$$\sum_{k_1+\ldots+k_d = n \atop{k_1,\ldots,k_d \in \{0,\ldots, n\}}} (k_1! k_2! \cdots k_d!)^{-2}.$$
To this end we consider the $d^{\text{th}}$ power of the series, i.e.,
$$\left(1 + \frac{x^{2\cdot 1}}{(1!)^2} + \frac{x^{2\cdot 2}}{(2!)^2} + \frac{x^{2\cdot 3}}{(3!)^2} + \ldots \right)^d.$$
To understand how to expand this expression imagine the $d$ brackets written out as a product. We have to pick exactly one factor from each of the $d$ brackets, multiply these $d$ factors, and sum over all possible choices. Thereby each of the chosen factors has the form $\frac{x^{2k}}{(k!)^2}$, and if we multiply $d$ such factors, say $\frac{x^{2k_1}}{(k_1!)^2}, \frac{x^{2k_2}}{(k_2!)^2}, \ldots \frac{x^{2k_d}}{(k_d!)^2}$, we get
$$\frac{x^{2(k_1+k_2+\ldots+k_d)}}{(k_1!k_2!\cdots k_d!)^2}.$$
Now we see that the exponent of $x$ equals $2n$ if and only if $k_1+k_2+\ldots + k_d = n$, so we have to sum over all choices of $k_1,\ldots k_d \in \{0,\ldots,n\}$ satisfying this condition and obtain that the coefficient of $x^{2n}$ is given by
$$\sum_{k_1+\ldots+k_d = n \atop{k_1,\ldots,k_d \in \{0,\ldots, n\}}} \frac{1}{(k_1!k_2!\cdots k_d!)^2},$$
as desired.

\section{Concluding remarks and pitfalls}

We have assumed that all packs of Skittles contain the same number $n$ of candies, which is actually not the case. In \cite{PossiblyWrong} it is pointed out that, assuming the number $n$ of Skittles in a pack is independently distributed with probability density function $f$, the probability that two randomly purchased packs are identical is given by
$$\sum_{n=1}^{\infty} f(n)^2\, \mathbb{P}[E_d^n].$$
Moreover, it says that they guessed $f(n)$ based on similar past studies and thereby obtained an expected value of 400-500 packs until the first match, depending on the assumptions for the density $f$.

Concluding this article we want to point out a possible pitfall. It is relatively easy to determine the number of different packs of Skittles if each pack contains $n$ Skittles and $d$ colors are available, and one obtains $\binom{n+d-1}{d-1}$ possibilities.\footnote{Note that $\binom{n+d-1}{d-1}$ is the cardinality $m$ of the hyperplane $H$ considered in Section~\ref{sec:expected_number}.} For instance, if $n=60$ and $d=5$ this gives
$$\binom{64}{4} = 635,376$$
different packs. Concluding from this that the probability for two randomly purchased packs of Skittles to be identical equals $1/635,376$ would be wrong, because the different packs are not equally likely. For instance, a pack with only red Skittles is less likely than a pack with twelve Skittles of each color.

\section{Behind the curtain -- the formal set-up}\label{sec:set up}

In this article we computed probabilities in an intuitive manner, and we want to specify this here by stating the corresponding probability spaces precisely. For the $2$-dimensional random walks with $n$ steps considered in section \ref{Two colours}, which represent the filling processes of an $n$-pack of Skittles with only two possible colors, we use 
$$
 \Omega = \big\{(x_1, \ldots, x_n) \,:\, x_1,\ldots, x_n \in \{r,u\} \big\},
$$
where $x_i=r$ means that step $i$ is a step to the right, and $x_i=u$ indicates that the $i$th step is a step up. Then, we have $|\Omega|=2^n$ and $\mathbb{P}[A] = |A|/2^n$ for each $A \subseteq \Omega$, as used intuitively above.

For the $3$-dimensional case corresponding to $d=3$ colours we use $\Omega=\{(x_1,\ldots,x_n)\,:\, x_1,\ldots, x_n \in \{r,f,u\}\}$ with the obvious meanings of $x_i=r,f,u$, and $\mathbb{P}[A] = |A|/3^n$ for each $A \subseteq \Omega$, since $|\Omega|=3^n$. For the general case, we let
$$
  \Omega=\big\{(x_1,\ldots,x_n)~|~ x_1,\ldots, x_n \in \{1,2,\ldots, d\}\big\},
$$
where $x_i=k$ means that the $i$th step is a step in direction $x_k$ for each $k \in \{1,\ldots, d\}$. Clearly, we then have $|\Omega|=d^n$.

When asking ourselves with which probability two randomly and independently chosen $n$-packs of Skittles are identical (i.e., two $n$-step random walks have the same endpoint) we formally consider the sample space $\Omega^2$ of all pairs $(W_1,W_2)$ of random walks $W_1, W_2$. The probability measure is then the product measure given by $\mathbb{P}[A_1 \times A_2] = |A_1|\cdot |A_2|/|\Omega|^2$. In the general case ($d$ possible colours/$d$-dimensional space) this means $\mathbb{P}[E_n^d] = |E_n^d|/d^{2n}$, as we have already used intuitively above.

When considering sequences of random walks as in Section \ref{sec:expected_number}, we formally deal with the product sample space $\Omega^{\mathbb{N}}$. The random variable $Z$ considered there is precisely defined by
\[
Z: \Omega^{\mathbb{N}} \to \mathbb{N}, \qquad (W_i)_{i \in \mathbb{N}} \mapsto \min\big\{n\in\mathbb{N}~:~ W_n \sim W_i \textnormal{ for some } i<n\big\}.
\]
To estimate $\mathbb{P}[Z=\ell]$, for fixed $n,d\in\mathbb{N}$ and $\ell \ge 2$, we consider the events $A_{ij}=\{(W_k)_{k\in \mathbb{N}}~:~ W_i \sim W_j\}$ and let $p:= \mathbb{P}[A_{ij}] = \mathbb{P}[E_d^n]$ for $i \ne j$. Then, we have
$$\mathbb{P}[Z=\ell]= \mathbb{P}\Bigg[ \bigcap_{i,j<\ell\atop{i\ne j}} A_{ij}^c \Bigg] \cdot \mathbb{P} \Bigg[\bigcup_{i=1}^{\ell-1} A_{\ell i} ~\Bigg|~ \bigcap_{i,j<\ell\atop{i\ne j}} A_{ij}^c\Bigg],$$
where $A_{ij}^c$ denotes the complement of $A_{ij}$. Here we clearly see that independence of the events $A_{ij}$ would be necessary to work as in the first presented plausible, yet incorrect approach.

\subsection*{Acknowledgment}
Joscha Prochno is supported by the Austrian Science Fund (FWF) with the Project P32405 ``Asymptotic Geometric Analysis and Applications''.
We thank the anonymous referee for helpful suggestions that improved the presentation of this paper. We thank Michael's brother-in-law Friedrich Delgado for pointing out the internet blog \cite{PossiblyWrong} to him. We also thank Gunther Leobacher (Graz) for reading a preliminary version of this article and his helpful comments and suggestions.

\begin{flushright}

Joscha Prochno\\
Faculty of Computer Science and Mathematics\\
University of Passau\\
Innstrasse 33, 94032 Passau, Germany\\
joscha.prochno@uni-passau.de\\[4mm]

Michael Schmitz\\
University of Flensburg\\
Auf dem Campus 1, 24943 Flensburg, Germany\\
michael.schmitz@uni-flensburg.de
\end{flushright}


\begin{thebibliography}{1}
	
	\bibitem{Bloom}  D. M. Bloom. \textit{A Birthday Problem}. American Mathematical Monthly, 80, 1141-2. (1973)
	\bibitem{Feller} W. Feller. An Introduction to Probability Theory and Its Applications (Vol. 1, 3rd ed.), New York: John Wiley (1968).
	\bibitem{mathoverflow} Mathoverflow. Birthday inequality for non-uniform distributions for fixed collision probability. \url{https://mathoverflow.net/questions/257027/birthday-inequality-for-non-uniform-distributions-for-} \url{fixed-collision-probabilit}. Last access on 06/12/2021.
	\bibitem{Munford} A. G. Munford. \textit{A Note on the Uniformity Assumption in the Birthday Problem}. The American Statistician, Vol. 31, No. 3 (Aug., 1977), p. 119
	\bibitem{Nunnikhoven} T.. S. Nunnikhoven. \textit{A Birthday Problem Solution for Nonuniform Birth Frequencies}.  The American Statistician, Vol. 46, No. 4 (Nov., 1992), pp. 270-274
	\bibitem{Perkins} W. Perkins. \textit{Birthday Inequalities, Repulsion, and hard Spheres}. Proc. Amer. Math. Soc. 144 (2016), 2635-2649
	\bibitem{PossiblyWrong} PossiblyWrong. \textit{Follow-up: I found two identical packs of Skittles, among 468 packs with a total of 27,740 Skittles}.
	\texttt{https://possiblywrong.wordpress.com/2019/04/06/follow-up-i\\-found-two-identical-packs-of-skittles-among-468-packs-with\\-a-total-of-27740-skittles/}. Last access on 06/12/2021.
	\bibitem{Schwarz} W. Schwarz. \textit{Approximating the Birthday Problem}. The American Statistician, Vol. 42, No. 3 (Aug., 1988), pp. 195-196
	\bibitem{Wiener} M. J. Wiener. \textit{Bounds on Birthday Attack Times}. IACR Eprint archive, \url{http://eprint.iacr.org/2005/318} (2005).
	
\end{thebibliography}
\end{document}